\documentclass[journal]{IEEEtran}

\usepackage{mathtools}
\usepackage{amssymb,dsfont}

\usepackage[compress]{cite}
\RequirePackage[colorlinks,linkcolor=blue,citecolor=blue,urlcolor=blue]{hyperref}

\usepackage{subcaption}

\usepackage[standard]{ntheorem}
\newtheorem{claim}{Claim}

\newcommand\ALPHABET{\mathds}
\newcommand\VEC{\mathbf}
\newcommand\PR{\mathds{P}}
\newcommand\EXP{\mathds{E}\strut}
\newcommand\reals{\mathds{R}}
\newcommand\DEFINED{\coloneqq}
\newcommand\BLANK{\mathfrak {E}}

\begin{document}

\title{Sufficient conditions for the value function and optimal strategy
to be even and quasi-convex}
\author{Jhelum Chakravorty and Aditya Mahajan%
\thanks{The authors are with Department of Electrical Engineering, McGill University, Montreal, QC, H3A 064, Canada (email: \texttt{jhelum.chakravorty@mail.mcgill.ca}, \texttt{aditya.mahajan@mcgill.ca}). 

This research was funded through NSERC Discovery Accelerator Grant 493011. 
}}
\maketitle

\begin{abstract}
  Sufficient conditions are identified under which the value function and the
  optimal strategy of a Markov decision process (MDP) are even and
  quasi-convex in the state. The key idea behind these conditions is the
  following. First, sufficient conditions for the value function and optimal
  strategy to be even are identified. Next, it is shown that if the value
  function and optimal strategy are even, then one can construct a
  ``folded MDP'' defined only on the non-negative values of the state space.
  Then, the standard sufficient conditions for the value function and optimal
  strategy to be monotone are ``unfolded'' to identify sufficient conditions
  for the value function and the optimal strategy to be quasi-convex.
  The results are illustrated by using an example of power allocation in
  remote estimation.
\end{abstract}


\section{Introduction}

\subsection{Motivation}

Markov decision theory is often used to identify structural or qualitative
properties of optimal strategies. Examples include control limit strategies
in machine maintenance~\cite{Derman1963, Kolesar1966}, threshold-based
strategies for executing call options~\cite{Taylor1967,Merton1973}, 
and monotone strategies in queueing
systems~\cite{Sobel1974,Stidham1989}. In all of these models, the optimal
strategy is \emph{monotone} in the state, i.e., if $x > y$ then the action
chosen at $x$ is greater (or less) than or equal to the action chosen at
$y$. Motivated by this, general conditions under which the optimal strategy
is monotone in scalar-valued states are identified in~\cite{Serfozo1976,
Puterman:1994, White1980, Ross:1983,Heyman:1984, Stokey1989}. Similar conditions for vector-valued
states are identified in~\cite{Topkis:1978, Topkis:1998, Papadaki2007}.
General conditions under which the value function is increasing and convex
are established in~\cite{Smith2002}.

Most of the above results are motivated by queueing models where the state is the queue length which 
takes non-negative values. However, for typical
applications in systems and control, the state takes both positive and
negative values. Often, the system behavior is symmetric for positive and
negative values, so one expects the optimal strategy to be even. Thus, for
such systems, a natural counterpart of monotone functions are even and
quasi-convex (or quasi-concave) functions. In this paper, we identify
sufficient conditions under which the value function and optimal strategy are
even and quasi-convex. 

As a motivating example, consider a remote estimation system in which a
sensor observes a Markov process and decides whether to transmit the current
state of the Markov process to a remote estimator. There is a cost or
constraint associated with transmission. When the transmitter does not
transmit or when the transmitted packet is dropped due to interference, the
estimator generates an estimate of the state of the Markov process based on
the previously received states. The objective is to choose transmission and
estimation strategies that minimize either the expected distortion and cost
of communication or minimize expected distortion under the transmission
constraint. Variations of such models have been considered in~\cite{ImerBasar,LipsaMartins:2011,MH2012,NayyarBasarTeneketzisVeeravalli:2013,JC_AM_TAC17,JC-AM-ifac16,JC-AM17_arxiv}. 

In such models, the optimal transmission and estimation strategies are
identified in two steps. In the first step, the joint optimization of
transmission and estimation strategies is investigated and it is established
that there is no loss of optimality in restricting attention to estimation
strategies of a specific form. In the second step, estimation strategies are
restricted to the form identified in the first step and the structure of the
best response transmission strategies is established. In particular, it is
shown that the optimal transmission strategies are even and
quasi-convex.\footnote{When the action space is binary---as is the case in
  most of the models of remote estimation---an even and quasi-convex strategy
  is equivalent to one in which the action zero is chosen whenever the
absolute value of the state is less than a threshold; otherwise, action one is chosen.} Currently, in the
literature these results are established on a case by case basis. For
example, see~\cite[Theorem~1]{LipsaMartins:2011}, ~\cite[Theorem~3]{NayyarBasarTeneketzisVeeravalli:2013},~\cite[Theorem~1]{XuHes2004a},~\cite[Theorem~1]{JC_AM_TAC17} among others.

In this paper, we identify sufficient conditions for the value functions and
optimal strategy of a Markov decision process to be even and quasi-convex. We
then consider a general model of remote estimation and verify these sufficient
conditions.

\subsection{Model and problem formulation}

Consider a Markov decision process (MDP) with state space $\ALPHABET X$ 
(which is either $\reals$, the real line, or a symmetric subset of the form $[-a, a]$) and
action space $\ALPHABET U$ (which is either a countable set or a compact
subset of reals).

Let $X_t \in \ALPHABET X$ and $U_t \in \ALPHABET U$ denote the state and
action at time~$t$. The initial state $X_1$ is distributed according to the
probability density function $\mu$ and the state evolves in a controlled
Markov manner, i.e., for any Borel measurable subset $A$ of $\ALPHABET X$,
\begin{multline*}
  \PR(X_{t+1} \in A \mid X_{1:t} = x_{1:t}, U_{1:t} = u_{1:t})
  \\  = 
  \PR(X_{t+1} \in A \mid X_{t} = x_{t}, U_{t} = u_{t}),
\end{multline*}
where $x_{1:t}$ is a short hand notation for $(x_1, \dots, x_t)$ and a
similar interpretation holds for $u_{1:t}$. We assume that there exists a
(time-homogeneous) controlled transition density $p(y|x;u)$ which is
continuous in $y$ for any $x \in \ALPHABET X$ and $u \in \ALPHABET U$ and for
any Borel measurable subset $A$ of $\ALPHABET X$, 
\[
  \PR(X_{t+1} \in A \mid X_t = x, U_t = u) = 
  \int_{A} p(y|x;u) dy.
\]
We use $p(u)$ to denote transition density corresponding to action ${u \in \ALPHABET
U}$.

The system operates for a finite horizon~$T$. For any time $t \in \{1, \dots,
{T-1}\}$, a measurable and lower semi-continuous\footnote{A function is lower
semi-continuous if and only if its lower level sets are closed.} function $c_t \colon \ALPHABET X \times \ALPHABET U \to
\reals$ denotes the instantaneous cost at time~$t$ and at the terminal
time~$T$ a measurable and lower semi-continuous function $c_T \colon \ALPHABET X \to \reals$
denotes the terminal~cost. 

The actions at time~$t$ are chosen according to a Markov strategy~$g_t$,
i.e.,
\[
  U_t = g_t(X_t), \quad t \in \{1, \dots, {T-1} \}.
\]
The objective is to choose a decision strategy $\VEC g \DEFINED (g_1, \dots,
g_{T-1})$ to minimize the expected total cost
\begin{equation} 
  J_T(\VEC g) \DEFINED \EXP^{g}\Big[ 
    \sum_{t=1}^{T-1} c_t(X_t, U_t) + c_T(X_T) 
  \Big].
  \label{eq:cost}
\end{equation}
We denote such an MDP by $(\ALPHABET X, \ALPHABET U, p, c_t)$.

From Markov decision theory~\cite{HernandezLermaLasserre:1996}, we know that
an optimal strategy is given by the solution of the following dynamic
program. Recursively define value functions $V_t \colon \ALPHABET X \to
\reals$ and value-action functions $Q_t \colon \ALPHABET X \times \ALPHABET U
\to \reals$ as follows: for all $x \in \ALPHABET X$ and $u \in \ALPHABET U$, 
\begin{align}
  V_T(x) &= c_T(x) ,
  \label{eq:V_T}
  \\
  \intertext{and for $t \in \{ {T-1}, \dots, 1 \}$,}
  Q_t(x, u) &= c_t(x,u) + \EXP[ V_{t+1}(X_{t+1}) \mid X_t = x, U_t = u ] 
        \notag \\
        &= c_t(x,u) + \int_{\ALPHABET X} p(y|x;u)V_{t+1}(y) dy, 
  \label{eq:Q_t} \\
  V_t(x) &= \min_{u \in \ALPHABET U} Q_t(x,u).
  \label{eq:V_t}
\end{align}

Then, a strategy $\VEC g^* = (g^*_1, \dots, g^*_{T-1})$ defined as
\[
  g^*_t(x) \in \arg\min_{u \in \ALPHABET U} Q_t(x,u)
\]
is optimal. To avoid ambiguity when the arg min is not unique, we pick 
\begin{equation}\label{eq:g_star}
    g_t^*(x) =  
    \begin{dcases}
      \max \big\{ v \in \arg \min_{u \in \ALPHABET U} Q_t(x,u) \big\}, &
      \hbox{if $x \ge 0$} \\
      \min \big\{ v \in \arg \min_{u \in \ALPHABET U} Q_t(x,u) \big\}, &
      \hbox{if $x < 0$}.
    \end{dcases}
\end{equation}

Let $\ALPHABET X_{\ge 0}$ and $\ALPHABET X_{> 0}$ denote the sets $\{x \in
\ALPHABET X : x \ge 0\}$ and $\{x \in \ALPHABET X : x > 0 \}$. We say that a
function $f \colon \ALPHABET X \to \reals$ is \emph{even and quasi-convex} if
it is even and for $x,x' \in \ALPHABET X_{\ge 0}$ such that $x < x'$, we
have that $f(x) \le f(x')$. The main contribution of this paper is to
identify sufficient conditions under which $V_t$ and $g^*_t$ are even and
quasi-convex.

\subsection{Main result}

\begin{definition}
For a given $u \in \ALPHABET U$, we say that a controlled transition density
$p(u)$ on $\ALPHABET X \times \ALPHABET X$ is \emph{even} if for all $x, y \in
\ALPHABET X$, $p(y|x;u) = p(-y|{-x};u)$.
\end{definition}

Our main result is the following.
\begin{theorem}\label{thm:main}
  Given an MDP $(\ALPHABET X, \ALPHABET U, p, c_t)$, define for $x,y \in
  \ALPHABET X_{\ge 0}$ and $u \in \ALPHABET U$, 
  \begin{equation}\label{eq:S}
   S(y|x;u) \DEFINED 1 - \int_{A_y} [ p(z| x;u) + p(-z| x;u)] dz,
  \end{equation}
  where $A_y = \{z \in \ALPHABET X : z < y \}$.
  Consider the following conditions:
  \begin{enumerate}
    \item[\textup{(C1)}] $c_T(\cdot)$ is even and quasi-convex and for $t \in
      \{1,\dots,{T-1}\}$ and $u \in \ALPHABET U$, $c_t(\cdot, u)$ is even and
      quasi-convex.
    \item [\textup{(C2)}] For all $u \in \ALPHABET U$, $p(u)$ is even.
    \item [\textup{(C3)}] For all $u \in \ALPHABET U$ and $y \in \ALPHABET
      X_{\ge 0}$, $S(y|x;u)$ is increasing for $x \in \ALPHABET X_{\ge 0}$.
    \item [\textup{(C4)}] For $t \in \{1, \dots, T-1\}$, $c_t(x,u)$ is
      submodular\footnote{Submodularity is defined in
      Sec.~\ref{sec:strategy}.} in $(x,u)$ on $\ALPHABET X_{\ge 0} \times
        \ALPHABET U$. 
    \item [\textup{(C5)}] For all $y \in \ALPHABET X_{\ge 0}$, $S(y|x;u)$ is
      submodular
      in $(x,u)$ on $\ALPHABET X_{\ge 0} \times \ALPHABET U$.
  \end{enumerate}
  Then, under \textup{(C1)--(C3)}, $V_t(\cdot)$ is even and quasi-convex for
  all $t \in \{1,\dots,T\}$ and under \textup{(C1)--(C5)}, $g^*_t(\cdot)$ is
  even and quasi-convex for all $t \in \{1, \dots, T-1\}$. 
\end{theorem}

The main idea of the proof is as follows. First, we identify conditions under
which the value function and optimal strategy of an MDP are even. Next, we
show that if we construct an MDP by ``folding'' the transition density, then
the ``folded MDP'' has the same value function and optimal strategy as the
original MDP for non-negative values of the state. Finally, we show that
if we take the sufficient conditions under which the value function and the
optimal strategy of the folded MDP are increasing and ``unfold'' these
conditions back to the original model, we get conditions (C1)--(C5) above.
The details are given in Sections~\ref{sec:evenMDP} and~\ref{sec:monotonocity}.

\section{Even MDPs and folded representations}\label{sec:evenMDP}

We say that \emph{an MDP is even} if for every~$t$ and every $u \in \ALPHABET
U$, $V_t(x)$, $Q_t(x,u)$ and $g^*_t(x)$ are even in~$x$. We start by
identifying sufficient conditions for an MDP to be even.

\subsection{Sufficient conditions for MDP to be even}

\begin{proposition}\label{prop:V-even}
  Suppose an MDP $(\ALPHABET X, \ALPHABET U, p, c_t)$ satisfies the following properties:
  \begin{description}
    \item [\textup{(A1)}] $c_T(\cdot)$ is even and for every $t \in \{1,\cdots, T-1\}$ and $u \in \ALPHABET U$,
      $c_t(\cdot, u)$ is even.
    \item [\textup{(A2)}] For every $u \in \ALPHABET U$, the
      transition density $p(u)$ is even.
      \end{description}
  Then, the MDP is even. 
\end{proposition}

\begin{proof}
  We proceed by backward induction. $V_T(x) = c_T(x)$ which is even by (A1).
  This forms the basis of induction. Now assume that $V_{t+1}(x)$ is even in
  $x$. For any $u \in \ALPHABET U$, we show that $Q_t(x,u)$  is even in $x$.
  Consider,
  \begin{align*}
    &\hskip -2em Q_t(-x, u) = c_t(-x, u) + \int_{\ALPHABET X} p(y| {-x};u)V_{t+1}(y) dy\\
    &\hskip -2em \stackrel{(a)}=
    c_t(x,u) + \int_{\ALPHABET X} p({-z}| {-x};u) V_{t+1}(-z) dz \\
    &\hskip -2em \stackrel{(b)}=
    c_t(x,u) + \int_{\ALPHABET X} p(z| x;u) V_{t+1}(z) dz 
    = Q_t(x,u)
  \end{align*}
  where $(a)$ follows from (A1), a change of variables $y = -z$, and the
  fact that $\ALPHABET X$ is a symmetric interval; and
  $(b)$ follows from (A2) and the induction hypothesis that
  $V_{t+1}(\cdot)$ is even. Hence, $Q_t(\cdot, u)$ is even.

  Since $Q_t(\cdot, u)$ is even, Eqs.~\eqref{eq:V_t} and~\eqref{eq:g_star}
  imply that  $V_t$ and $g^*_t$ are also even. Thus, the result is true for
  time~$t$ and, by induction, true for all time~$t$.
\end{proof}

\subsection{Folding operator for distributions}
We now show that if the value function is even, we can construct a ``folded''
MDP with state-space $\ALPHABET X_{\ge 0}$ such that the value function and
optimal strategy of the folded MDP match that of the original MDP on
$\ALPHABET X_{\ge 0}$. For that matter, we first define the following:

\begin{definition}[Folding Operator] \label{def:folding_prob}
Given a probability density $\pi$ on $\ALPHABET X$, the folding operator
$\mathcal F \pi$ gives a density $\tilde \pi$ on $\ALPHABET X_{\ge 0}$ such
that for any $x \in \ALPHABET X_{\ge 0}$, $\tilde \pi(x) = \pi(x) + \pi(-x)$.
\end{definition}

As an immediate implication, we have the following:
\begin{lemma}\label{lem:even-folding}
  If $f: \ALPHABET X \to \reals$ is even, then for any probability
  distribution $\pi$ on $\ALPHABET X$ and $\tilde \pi = \mathcal F \pi$, we
  have 
  \[
    \int_{\ALPHABET X} f(x) \pi(x) dx
    = \int_{\ALPHABET X_{\ge 0}} f(x) \tilde \pi(x) dx.
  \]
\end{lemma}

Now, we generalize the folding operator to transition densities.

\begin{definition}\label{def:folding_mat}
  Given a transition density $p$ on $\ALPHABET X \times \ALPHABET
  X$, the folding operator $\mathcal F p$ gives a transition
  density $\tilde p$ on $\ALPHABET X_{\ge 0} \times \ALPHABET X_{\ge 0}$ such
  that for any $x,y \in \ALPHABET X_{\ge 0}$, 
  \(
    \tilde p(y | x) = p(y|x) + p(-y|x).
  \)
\end{definition}

\begin{definition}[Folded MDP]\label{def:folded}
  Given an MDP $(\ALPHABET X, \ALPHABET U, p, c_t)$, define the \emph{folded
  MDP} as $(\ALPHABET X_{\ge 0}, \ALPHABET U, \tilde p, c_t)$, where for all
  $u \in \ALPHABET U$, $\tilde p(u) = \mathcal F p(u)$.
\end{definition}

Let $\tilde Q_t$ and $\tilde V_t$ and $\tilde g^*_t$ denote respectively the
value-action function, the value function, and the optimal strategy of the
folded MDP. Then, we have the following.

\begin{proposition}\label{prop:folded}
  If the MDP $(\ALPHABET X, \ALPHABET U, p, c_t)$ is even, then for any $x
  \in \ALPHABET X$ and $u \in \ALPHABET U$, 
  \begin{equation}\label{eq:tilde-Q-V-g}
    Q_t(x,u) = \tilde Q_t(|x|,u),
    \quad 
    V_t(x) = \tilde V_t(|x|),
    \quad 
    g^*_t(x) = \tilde g^*_t(|x|).
  \end{equation}
\end{proposition}

\begin{proof}
  We proceed by backward induction. For $x \in \ALPHABET X$ and $\tilde x \in
  \ALPHABET X_{\ge 0}$, $V_T(x) = c_T(x)$ and $\tilde V_T(\tilde x) =
  c_T(\tilde x)$. Since $V_T(\cdot)$ is even, $V_T(x) = V_T(|x|) = \tilde
  V_T(|x|)$. This is the basis of induction. Now assume that for all $x \in
  \ALPHABET X$, $V_{t+1}(x) = \tilde V_{t+1}(|x|)$. Consider $x \in \ALPHABET
  X_{\ge 0}$ and $u \in \ALPHABET U$. Then we have
  \begin{align*}
   &\hskip -2em  Q_t(x,u) = c_t(x,u) + \int_{\ALPHABET X} p(y| x;u) V_{t+1}(y) dy\\
    &\hskip -2em \stackrel{(a)} = c_t(x,u) + \int_{\ALPHABET X_{\ge 0}} \tilde p(y| x;u) V_{t+1}(y) dy\\
    & \hskip -2em \stackrel{(b)} = c_t(x,u) + \int_{\ALPHABET X_{\ge 0}} \tilde p(y| x;u) \tilde V_{t+1}(y) dy
    = \tilde Q_t(x,u),
  \end{align*}
  where $(a)$ uses Lemma~\ref{lem:even-folding} and that $V_{t+1}$ is even
  and $(b)$ uses the induction hypothesis. 

  Since the $Q$-functions match for $x \in \ALPHABET X_{\ge 0}$,
  equations~\eqref{eq:V_t} and~\eqref{eq:g_star} imply that the value
  functions and the optimal strategies also match on $\ALPHABET X_{\ge 0}$,
  i.e., for $x \in \ALPHABET X_{\ge 0}$, 
  \[
    V_t(x) = \tilde V_t(x) 
    \quad \text{and} \quad 
    g^*_t(x) = \tilde g^*_t(x).
  \]
  Since $V_t$ and $g^*_t$ are even, we get that~\eqref{eq:tilde-Q-V-g} is
  true at time $t$. Hence, by principle of induction, it is true for all $t$.
\end{proof}

\section{Monotonicity of the value function and the optimal strategy}\label{sec:monotonocity}

We have shown that under (A1) and (A2) the original MDP is equivalent to a
folded MDP with state-space $\ALPHABET X_{\ge 0}$. Thus, we can use standard
conditions to determine when the value function and the optimal strategy of
the folded MDP are monotone. Translating these conditions back to the
original model, we get the sufficient conditions for the original model.

\subsection{Monotonicity of the value function}

The results on monotonicity of value functions rely on the notion of
stochastic monotonicity.

Given a transition density $p$ defined on $\ALPHABET X$, the cumulative
transition distribution function $P$ is defined as
  \[
    P(y|x) = \int_{A_y} p(z|x) dz, 
    \quad \hbox{where }
    A_y = \{z \in \ALPHABET X : z < y \}.
  \]

\begin{definition}[Stochastic Monotonicity] \label{def:stoc-mono}
  A transition density $p$ on $\ALPHABET X$ is said to be
  \emph{stochastically monotone increasing} if for every $y \in \ALPHABET X$, the
  cumulative distribution function $P(y|x)$ corresponding to $p$ is decreasing
  in~$x$.
\end{definition}

\begin{proposition} \label{prop:value}
  Suppose the folded MDP $(\ALPHABET X_{\ge 0}, \ALPHABET U, \tilde p, c_t)$
  satisfies the following:
  \begin{description}
    \item [\textup{(B1)}] $c_T(x)$ is increasing in $x$ for $x \in \ALPHABET X_{\ge
      0}$; for any $t \in \{1, \dots, T - 1\}$ and $u \in \ALPHABET U$,
      $c_t(x,u)$ is increasing in $x$ for $x \in \ALPHABET X_{\ge 0}$. 
    \item [\textup{(B2)}] For any $u \in \ALPHABET U$, $\tilde p(u)$ is
      stochastically monotone increasing.
  \end{description}
  Then, for any $t \in \{1, \dots, T\}$, $\tilde V_t(x)$ is increasing in $x$ for $x
  \in \ALPHABET X_{\ge 0}$. 
\end{proposition}
A version of this proposition when $\ALPHABET X$ is a subset of integers is
given in~\cite[Theorem~4.7.3]{Puterman:1994}. The same proof argument also
works when $\ALPHABET X$ is a subset of reals.

Recall the definition of $S$ given in~\eqref{eq:S}. 
(B2) is equivalent to the following: 
\begin{description}
  \item[(B2')] For every $u \in \ALPHABET U$ and $x,y \in \ALPHABET X_{\ge
    0}$, $S(y|x;u)$ is increasing in $x$.
\end{description}

An immediate consequence of Propositions~\ref{prop:V-even},
\ref{prop:folded}, and~\ref{prop:value} is the following:
\begin{corollary} \label{cor:value}
  Under \textup{(A1)}, \textup{(A2)}, \textup{(B1)}, and \textup{(B2)} (or
  \textup{(B2')}), the value functions $V_t(\cdot)$ is even
  and quasi-convex. 
\end{corollary}

\begin{remark}\label{rem:1}
  Note that (A1) and (B1) are equivalent to (C1), (A2) is same as (C2),
  and (B2) (or equivalently, (B2')) is equivalent to (C3). Thus, Corollary~\ref{cor:value}
  proves the first part of Theorem~\ref{thm:main}.
\end{remark}

\subsection{Monotonicity of the optimal strategy} \label{sec:strategy}
Now we state sufficient conditions under which the optimal strategy is
increasing. These results rely on the notion of submodularity.

\begin{definition}[Submodular function]
  A function $f \colon \ALPHABET X \times \ALPHABET U \to \reals$ is called 
  submodular if for any $x, y \in \ALPHABET X$ and $u, v \in \ALPHABET U$
  such that $x \ge y$ and $u \ge v$, we have
  \[
    f(x,u) + f(y,v) \le f(x, v) + f(y, u).
  \]
\end{definition}

An equivalent characterization of submodularity is that 
\begin{align*}
  f(y, u) - f(y, v) \ge f(x, u) - f(x, v),\\
  \implies f(x, v) - f(y, v) \ge f(x, u) - f(y, u),
\end{align*}
which implies that the differences in one variable are decreasing in the other.

\begin{proposition}\label{prop:strategy}
  Suppose that in addition to \textup{(B1)} and \textup{(B2)} (or
  \textup{(B2')}), the folded MDP $(\ALPHABET X_{\ge 0}, \ALPHABET U, \tilde
  p, c_t)$ satisfies the following property:
  \begin{description}
    \item [\textup{(B3)}] For all $t \in \{1, \dots, T-1\}$, $c_t(x,u)$ is submodular
      in $(x,u)$ on $\ALPHABET X_{\ge 0} \times \ALPHABET U$. 
    \item [\textup{(B4)}] For all $y \in \ALPHABET X_{\ge 0}$, $S(y|x;u)$ is submodular
      in $(x,u)$ on $\ALPHABET X_{\ge 0} \times \ALPHABET U$, where $S(y|x;u)$
      is defined in~\eqref{eq:S}.
  \end{description}
  Then, for every $t \in \{1,\cdots,T-1\}$, the optimal strategy $\tilde
  g^*_t(x)$  is increasing in $x$ for $x \in \ALPHABET X_{\ge 0}$.
\end{proposition}
A version of this proposition when $\ALPHABET X$ is a subset of integers is
given in~\cite[Theorem~4.7.4]{Puterman:1994}. The same proof argument also
works when $\ALPHABET X$ is a subset of reals.

An immediate consequence of Propositions~\ref{prop:V-even},
\ref{prop:folded}, \ref{prop:value}, and~\ref{prop:strategy} is the following:
\begin{corollary} \label{cor:strategy}
  Under \textup{(A1)}, \textup{(A2)}, \textup{(B1)}, \textup{(B2)} (or
  \textup{(B2')}), \textup{(B3)}, and \textup{(B4)} the optimal strategy
  $g^*_t(\cdot)$ is even and quasi-convex. 
\end{corollary}
\begin{remark} \label{rem:2}
  As argued in Remark~\ref{rem:1}, (A1), (A2), (B1), (B2) (or (B2')) are equivalent to
  (C1)--(C3). Note that (B3), (B4) is the same as (C4), (C5). Thus,
  Corollary~\ref{cor:strategy} proves the second part of
  Theorem~\ref{thm:main}.
\end{remark}

\section{Remark on infinite horizon setup}\label{sec:infin}

Although we restricted attention to finite horizon models, the results extend immediately to infinite horizon discounted cost setup. In particular, suppose the per-step cost is time-homogeneous and given by $c \colon \ALPHABET X \times \ALPHABET U \to \reals$ and future is discounted by $\beta \in (0,1)$. Define the following Bellman operators: for any $g \colon \ALPHABET X \to \ALPHABET U$, and $V \colon \ALPHABET X \to \reals$
\[
[{\mathcal B}_g V](x) = c(x,g(x)) + \beta \int_X p(y|x;g(x)) V(y) dy
\]
and 
\[
\mathcal B^* V = \min_{g \colon \ALPHABET X \to \ALPHABET U} \mathcal B_g V.
\]

Suppose the model satisfies standard conditions (see~\cite[Chapter~4]{HernandezLermaLasserre:1996}) so that $\mathcal B^*$ is a contraction and has a unique fixed point (which we denote by $V_\beta$) and there exists a strategy $g_\beta \colon \ALPHABET X \to \ALPHABET U$ such that $V_\beta = \mathcal B_{g_\beta} V_\beta$. Then, the result of Theorem~\ref{thm:main} is also true for $V_\beta$ and $g_\beta$. In particular,

\begin{corollary}
 Given an MDP $(\ALPHABET X, \ALPHABET U, p, c)$ and a discount factor $\beta \in (0,1)$, consider the following conditions:

  \begin{enumerate}
    \item[\textup{(C1')}] For $u \in \ALPHABET U$, $c(\cdot, u)$ is even and
      quasi-convex.
    \item [\textup{(C4')}] $c(x,u)$ is
      submodular in $(x,u)$ on $\ALPHABET X_{\ge 0} \times
        \ALPHABET U$. 
      submodular
      in $(x,u)$ on $\ALPHABET X_{\ge 0} \times \ALPHABET U$.
  \end{enumerate}
  Then, under \textup{(C1'), (C2), (C3)}, $V_\beta(\cdot)$ is even and quasi-convex and under \textup{(C1'), (C2), (C3), (C4'), (C5)}, $g^*_\beta(\cdot)$ is
  even and quasi-convex. 
\end{corollary}
\begin{proof}
Note that the equivalence to folded MDP continues to hold for infinite horizon setup. Therefore, the result follows from extension of Propositions 3 and 4 to infinite horizon setup. For example, see~\cite[Section~6.11]{Puterman:1994}.
\end{proof}

\section {Remarks about discrete $\ALPHABET X$}

So far we assumed that $\ALPHABET X$ was a subset of the real line. Now
suppose $\ALPHABET X$ is discrete (either the set $\mathbb Z$ of integers or
a symmetric subset of the form $\{-a, \dots, a\}$). With a slight abuse of
notation, let $p(y|x;u)$ denote $\PR(X_{t+1} = y | X_t = x, U_t = u)$. 

\begin{theorem} \label{thm:discrete}
  The result of Theorem~\ref{thm:main} is true for discrete $\ALPHABET X$
  with $S$ defined as
  \[
    S(y|x;u) = 1 - \sum_{z \in A_y}
    \big[ p(z|x;u) + p(-z|x;u) \big]
  \]
  where $A_y = \{x \in \ALPHABET X : x < y \}$.
\end{theorem}

The proof proceeds along the same lines as the proof of
Theorem~\ref{thm:main}. In particular, 
\begin{itemize}
  \item Proposition~\ref{prop:V-even} is also true for discrete $\ALPHABET X$.

  \item Given a probability mass function $\pi$ on
    $\ALPHABET X$, define the folding operator $\mathcal F$  as follows:
    $\tilde \pi = \mathcal F \pi$ means that $\tilde \pi(0) = \pi(0)$ and for
    any $x \in \ALPHABET X_{> 0}$, $\tilde \pi(x) = \pi(x) + \pi(-x)$. 

  \item Use this definition of the folding operator to define the folded MDP,
    as in Definition~\ref{def:folded}. Proposition~\ref{prop:folded} remains
    true with this modified definition.

  \item A discrete state Markov chain with transition function $p$ is
    stochastically monotone increasing if for every $y \in \ALPHABET X$, 
    \[
      P(y|x) = \sum_{z \in \ALPHABET A_y} p(z|x),
      \quad \hbox{where }
      A_y = \{z \in \ALPHABET X : z < y \},
    \]
    is decreasing in $x$.
  \item Propositions~\ref{prop:value} and~\ref{prop:strategy} are also true
    for discrete $\ALPHABET X$. 
  \item The result of Theorem~\ref{thm:discrete} follows from
    Corollaries~\ref{cor:value} and~\ref{cor:strategy}.
\end{itemize}

\subsection{Monotone dynamic programming} \label{sec:monotoneDP}

Under (C1)--(C5), the even and quasi-convex property of the optimal strategy can be used to simplify the dynamic program given by~\eqref{eq:V_T}--\eqref{eq:V_t}. For conciseness, assume
that the state space $\ALPHABET X$ is a set of integers form $\{-a, -a+1,
\cdots,a-1, a\}$ and the action space $\ALPHABET U$ is a set of integers of
the form $\{\underline u, \underline u+1,\cdots, \bar u -1, \bar u\}$. 

Initialize $V_T(x)$ as in~\eqref{eq:V_T}. Now, suppose $V_{t+1}(\cdot)$ has
been calculated. Instead of computing $Q_t(x,u)$ and $V_t(x)$ according
to (appropriately modified versions of) \eqref{eq:Q_t} and~\eqref{eq:V_t}, we
proceed as follows:

\begin{enumerate}
  \item Set $x = 0$ and $w_x = \underline u$. 
  \item For all $u \in [w_x, \bar u]$, compute $Q_t(x,u)$ according
    to~\eqref{eq:Q_t}. 
  \item Instead of~\eqref{eq:V_t}, compute
    \begin{align*}
      V_t(x) &= \min_{u \in [w_x, \bar u]} Q_t(x,u), \quad \text{and set}\\
      g_t(x) &= \max\{ v \in [w_x, \bar u] \text{ s.t. } V_t(x) = Q_t(x,v) \}. 
    \end{align*}
  \item Set $V_t(-x) = V_t(x)$ and $g_t(-x) = g_t(x)$.
  \item If $x = a$, then stop. Otherwise, set $w_{x+1} = g_t(x)$ and $x
    = x + 1$. Go to step 2. 
\end{enumerate}

\subsection{A remark on randomized actions}

Suppose $\ALPHABET U$ is a discrete set of the form $\{\underline u,
\underline u + 1, \dots, \bar u \}$. In constrained optimization problems, it
is often useful to consider the action space $\ALPHABET W = [\underline u,
\bar u]$, where for $u, {u+1} \in \ALPHABET U$, an action $w \in (u, u+1)$
corresponds to a randomization between the ``pure'' actions $u$ and $u+1$.
More precisely, let transition probability $\breve p$ corresponding to
$\ALPHABET W$ be given as follows: for any $x,y \in \ALPHABET X$ and $w \in
(u, u+1)$,
\[
  \breve p(y|x;w) = (1 - \theta(w)) p(y|x;u) + \theta(w) p(y|x;u+1)
\]
where $\theta : \ALPHABET W \to [0, 1]$ is such that
for any $u \in \ALPHABET U$, 
\begin{equation}\label{eq:theta}
\lim_{w \downarrow u} \theta(w) = 0, \quad \text{and} \quad
\lim_{w \uparrow u+1} \theta(w) = 1.
\end{equation} 
Thus, $\breve p(w)$ is continuous at
all $u \in \ALPHABET U$.  

\begin{theorem}
  If $p(u)$ satisfies \textup{(C2)}, \textup{(C3)}, and \textup{(C5)} then so
  does $\breve p(w)$. 
\end{theorem}

\begin{proof}
  Since $\breve p(w)$ is linear in $p(u)$ and $p(u+1)$, both of which satisfy
  (C2) and (C3), so does $\breve p(w)$. 

  To prove that $\breve p(w)$ satisfies (C5), note that 
  \[
    \breve S(y|x,w) = S(y|x;u) + \theta(w)[ S(y|x,u+1) - S(y|x;u) ].
  \]
  So, for $v,w \in (u, u+1)$ such that $v > w$, we have that
  \[
    \breve S(y|x;v) - \breve S(y|x;w) = 
    \big( \theta(v) - \theta(w) \big)
    [ S(y|x; u+1) - S(y|x;u) ]
  \]
  Since $\theta(\cdot)$ is increasing, $\theta(v) - \theta(w) \ge 0$.
  Moreover, since $S(y|x;u)$ is submodular in $(x,u)$, $S(y|x; u+1) - S(y|x;u)$ is decreasing in $x$, and, therefore, so is $\breve S(y|x;v) - \breve
  S(y|x;w)$. Hence, $\breve S(y|x;w)$ is submodular in
  $(x,w)$ on $\ALPHABET X \times (u, u+1)$. Due to~\eqref{eq:theta}, $\breve
  S(y|x;w)$ is continuous in $w$. Hence, $\breve S(y|x;w)$ is submodular in $(x,w)$ on $\ALPHABET
  X \times [u, u+1]$. By piecing intervals of the form $[u, u+1]$ together,
  we get that $\breve S(y|x;w)$ is submodular on $\ALPHABET X \times
  \ALPHABET W$.
\end{proof}

\section{An example: Optimal power allocation strategies in remote estimation}\label{sec:example-RE}

Consider a remote estimation system that consists of a sensor and an
estimator. The sensor observes a first order autoregressive process
$\{X_t\}_{t \ge 1}$, $X_t \in \ALPHABET X$, where $\ALPHABET X$ is either $\reals$ or $\mathbb Z$. The system starts 
with $X_1 = 0$ and for $t > 1$, 
\[
  X_{t+1} = a X_t + W_t,
\]
where $a \in \ALPHABET X$ is a constant and $\{W_t\}_{t \ge 1}$, $W_t \in \ALPHABET X$ is an i.i.d.\
noise process with probability mass/density function $\varphi$. 

At each time step, the sensor uses power $U_t$ to send a packet containing
$X_t$ to the remote estimator. $U_t$ takes values in $[0, u_{\max}]$, where
$U_t = 0$ denotes that no packet is sent. The packet is received with
probability $q(U_t)$, where $q$ is an increasing function with $q(0) = 0$ and
$q(u_{\max}) \le 1$.

Let $Y_t$ denote the received symbol. $Y_t = X_t$ if the packet is received
and $Y_t = \BLANK$ if the packet is not received. Packet reception is
acknowledged, so the sensor knows $Y_t$ with one unit delay. At each stage,
the receiver generates an estimate $\hat X_t$ as follows. $\hat X_0$ is $0$
and for $t > 0$, 
\begin{equation}\label{eq:AR_hatX}
  \hat X_t = \begin{cases}
    a \hat X_{t-1}, & \hbox{if $Y_t = \BLANK$} \\
    Y_t, & \hbox{if $Y_t \neq \BLANK$}.
  \end{cases}
\end{equation}
Under some conditions, such an estimation rule is known to be
optimal~\cite{LipsaMartins:2011,
NayyarBasarTeneketzisVeeravalli:2013,shi2012optimal,JC-AM-ifac16,RenJohanssonarxiv2016, JC-AM17_arxiv}\footnote{The model presented above appears as an intermediate step in the analysis of remote estimation problem. One typically starts with a model where the transmission strategy is of the form $U_t = g_t(X_{1:t}, Y_{1:t-1}, U_{1:t-1})$ and the estimation strategy is of the form $\hat X_t = h_t(Y_{1:t})$. This is a decentralized control problem. After a series of simplifications, it is shown that there is no loss of optimality  to restrict attention to estimation strategies of the form~\eqref{eq:AR_hatX} (see~\cite[Fact~B.3]{LipsaMartins:2011}, ~\cite[Theorem~3]{NayyarBasarTeneketzisVeeravalli:2013},~\cite[Theorem~1]{JC-AM17_arxiv} among others).
Once the attention is restricted to estimation strategies of the form~\eqref{eq:AR_hatX}, the next step is to simplify the structure of the optimal transmission strategy (see~\cite[Fact~A.4]{LipsaMartins:2011},~\cite[Theorem~3]{NayyarBasarTeneketzisVeeravalli:2013},~\cite[Theorem~1]{XuHes2004a},~\cite[Theorem~1]{JC-AM17_arxiv} among others).
The model presented above corresponds to this step.}. 

There are two types of costs: (i)~a communication cost $\lambda(U_t)$, where
$\lambda$ is an increasing function with $\lambda(0) = 0$; and (ii)~an
estimation cost $d(X_t - \hat X_t)$, where $d$ is an even and quasi-convex
function with $d(0) = 0$.  

Define the error process $\{E_t\}_{t \ge 0}$ as 
\(
  E_t = X_t - a \hat X_{t-1}.
\)
The error process $\{E_t\}_{t \ge 0}$ evolves in a controlled Markov manner
as follows:
\begin{equation} \label{eq:dynamics}
  E_{t+1} = \begin{cases}
    a E_t + W_t, & \hbox{if $Y_t = \BLANK$} \\
    W_t, & \hbox{if $Y_t \neq \BLANK$}
  \end{cases}
\end{equation}
Due to packet acknowledgments, $E_t$ is measurable at the sensor at
time~$t$. If a packet is received, then $\hat X_t = X_t$ and the estimation
cost is $0$. If the packet is dropped, $X_t - \hat X_t = E_t$ and an
estimation cost of $d(E_t)$ is incurred.

The objective is to choose a transmission strategy $g = (g_1,
\dots, g_T)$ of the form
\(
  U_t = g_t(E_t)
\)
to minimize 
\[
  \EXP\bigg[ \sum_{t=1}^T \big[ \lambda(U_t) + (1- q(U_t))d(E_t) \big]
  \bigg].
\]

The above model is Markov decision process with state $E_t \in \ALPHABET X$,
control action $U_t \in [0, u_{\max}]$, per-step cost 
\begin{equation} \label{eq:cost}
  c(e,u) = \lambda(u) + (1 - q(u)) d(e),
\end{equation}
and transition density/mass function
\begin{equation} \label{eq:transition}
  p(e_{+}|e;u) = 
  q(u) \varphi(e_{+})
  + (1 - q(u)) \varphi(e_{+} - ae).
\end{equation}

For ease of reference, we restate the assumptions imposed on the cost:
\begin{enumerate}
  \item[(M0)] $q(0) = 0$ and $q(u_{\max}) \le 1$.
  \item[(M1)] $\lambda(\cdot)$ is increasing with $\lambda(0) = 0$.
  \item[(M2)] $q(\cdot)$ is increasing.
  \item[(M3)] $d(\cdot)$ is even and quasi-convex with $d(0) = 0$.
\end{enumerate}
In addition, we impose the following assumptions on the probability
density/mass function of the i.i.d.\
process $\{W_t\}_{t \ge 1}$:
\begin{enumerate}
  \item[(M4)] $\varphi(\cdot)$ is even.
  \item[(M5)] $\varphi(\cdot)$ is unimodal (i.e., quasi-concave).
\end{enumerate}

\begin{claim} \label{clm:prop}
  We have the following:
  \begin{enumerate}
    \item under assumptions \textup{(M0)} and \textup{(M3)}, the
      per step cost function given by~\eqref{eq:cost} satisfies \textup{(C1)}.
    \item under assumptions \textup{(M0)}, \textup{(M2)} and \textup({M3)},
      the per step cost function given by~\eqref{eq:cost} satisfies
      \textup{(C4)}. \item under assumption \textup{(M4)}, the transition
      density $p(u)$ given
      by~\eqref{eq:transition} satisfies \textup{(C2)}.
    \item under assumptions \textup{(M0)}, \textup{(M2)}, \textup{(M4)} and
      \textup{(M5)}, the transition density $p(u)$ satisfies \textup{(C3)}
      and \textup{(C5)}.
  \end{enumerate}
\end{claim}
The proof is given in Appendix~\ref{app:proof_clm_prop}.

An immediate consequence of Theorem~\ref{thm:main} and Claim~\ref{clm:prop}
is the following:
\begin{theorem} \label{thm:re}
  Under assumptions \textup{(M0)}, \textup{(M2)}--\textup{(M5)}, the value
  function and the optimal strategy for the remote estimation model are even
  and quasi-convex.
\end{theorem}

\begin{remark}
  Although Theorem~\ref{thm:re} is derived for continuous action space, it is
  also true when the action space is a discrete set. In particular, if we
  take the action space to be $\{0,1\}$ and $q(1) = 1$, we get the results of~\cite[Theorem~1]{LipsaMartins:2011},~\cite[Proposition~1]{ImerBasar},~\cite[Theorem~3]{NayyarBasarTeneketzisVeeravalli:2013},~\cite[Theorem~1]{JC_AM_TAC17}; if we take the action space to be
  $\{0,1\}$ and $q(1) = \varepsilon$, we get the result of~\cite[Theorem~1]{JC-AM-ifac16},~\cite[Theorem~2]{JC-AM17_arxiv}.
\end{remark}

\begin{figure}
  \centering
  \includegraphics[width=\linewidth]{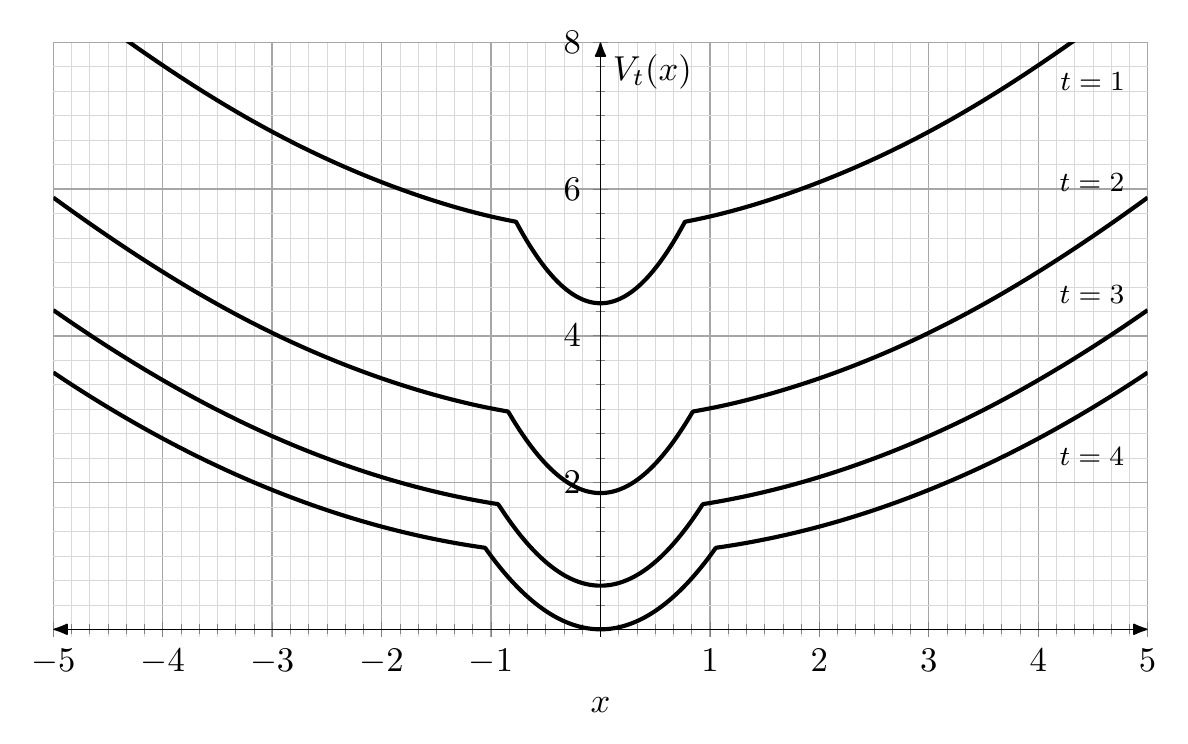}
  \caption{The value function for the remote estimation problem of horizon
    $T=4$ with state space $\ALPHABET X = \reals$, action space $\ALPHABET U =
    \{0, 1\}$, $a = 1$, $W_t \sim \mathcal{N}(0,1)$, $d(e) = e^2$, $\lambda =
    1$, $q(0) = 0$, and $q(1) = 0.9$. The kink in the value function
    corresponds to the point where the optimal action changes. The value
  functions of the folded MDP are identical to the value functions above when
restricted to the domain $\reals_{\ge 0}$.}
  \label{fig:value}
\end{figure}

To illustrate the above result, consider the case when $\ALPHABET X =
\reals$, $\ALPHABET U = \{0, 1\}$, $a = 1$, $W_t \sim \mathcal{N}(0,1)$, $d(e) = e^2$,
$\lambda = 1$, $q(0) = 0$, $q(1) = 0.9$, and $T = 4$. We discretize the state
space with a uniform grid of width $0.01$ and numerically solve the resulting
dynamic program~\eqref{eq:V_T}--\eqref{eq:V_t}. The value functions
across time are shown in Fig.~\ref{fig:value}. The optimal strategy is of the
form
\[
  g_t(e) = \begin{cases}
    1, & \hbox{if $|e| > k_t$} \\
    0, & \hbox{if $|e| \le k_t$}
  \end{cases}
\]
where $k_1 = 0.77$, $k_2 = 0.84$, $k_3 = 0.93$, and $k_4 = 1.05$. Note that,
as expected, both the value function and the optimal policy and even and
quasi-convex. 

\subsection{Some comments on the conditions}

Note that the result does not depend on (M1). This is for the
following reason. Suppose there are two power
levels $u_1, u_2 \in \ALPHABET U$ such that $u_1 < u_2$ but $\lambda(u_1)
\ge \lambda(u_2)$, then for any $e \in \ALPHABET X$, $c(e, u_1) \ge c(e,
u_2)$. Thus, action~$u_1$ is dominated by action~$u_2$ and is, therefore,
never optimal and can be eliminated.

All other conditions, (M0), (M2)--(M5) are also necessary as is explained
below. 

Condition (M2) is necessary. We illustrate that with the following example.
Consider an example where $\ALPHABET X = \mathbb Z$ and $\ALPHABET U = \{0,
u_1, u_2\}$ such that $u_1 < u_2$ but $q(u_1) > q(u_2)$. Then, we can
consider an alternative action space $\ALPHABET U' = \{0, u'_1, u'_2\}$ where
$u'_1 < u'_2$ and the bijection $\sigma : \ALPHABET U \to \ALPHABET U'$ such
that $\sigma(0) = 0$, $\sigma(u_1) = u'_2$ and $\sigma(u_2) = u'_1$. Now
consider a remote estimation system with communication cost $\lambda' =
\lambda \circ \sigma^{-1}$ and success probabilities $q' = q \circ
\sigma^{-1}$. By construction $q'$ satisfies (M0) and
(M2).\footnote{$\lambda'$ does not satisfy (M1), but (M1) is not needed for
Theorem~\ref{thm:re}.} If $d(\cdot)$ and $\varphi(\cdot)$ are chosen to
satisfy (M3)--(M5), then by Theorem~\ref{thm:re}, the optimal strategy $g'
\colon \ALPHABET X \to \ALPHABET U'$ is even and quasi-convex. In particular,
we can pick $\lambda$, $d$ and $\varphi$ such that $g'(0) = 0$, $g'(\pm1) =
u'_1$ and $g'(\pm 2) = u'_2$. However, this means that with the original
labels, the optimal strategy would have been $g = g' \circ \sigma^{-1}$, which
means $g(0) = 0$, $g(\pm 1) = u_2$ and $g(\pm 2) = u_1$. And hence, the
optimal strategy is not quasi-convex.

Conditions (M3) and (M4) are necessary. If they are not satisfied, then it is
easy to construct examples where the value function is not even. 

Condition (M5) is also necessary. We illustrate that with the following example.
Consider an example where $\ALPHABET X = \mathbb Z$. In particular,
let $a = 1$ and $\varphi$ have support $\{-1, 0, 1\}$ where $\varphi(0) =
1-2p$ and $\varphi(-1) = \varphi(1) = p$. Suppose $p > 1/3$, so that (M5) is
not satisfied. Furthermore, suppose $T = 2$, $\ALPHABET U = \{0, 1\}$ and
consider the following functions: $\lambda(0) = 0$,
$\lambda(1) = K$; $q(0) = 0$ and $q(1) = 1$; and $d(0) = 0$, $d(\pm 1) = 1$,
and for any $e \not\in \{-1, 0, 1\}$, $d(e) = 1 + k$, where $k$ is a positive
constant. Note that $q(\cdot)$ satisfies (M0) and (M2); $d(\cdot)$ satisfies
(M3); and $\varphi(\cdot)$ satisfies (M4) but not (M5). Suppose $K >  2(1 +
k)$, so that action $1$ is not optimal at any time. Thus, $V_2(e) = d(e)$ and
$V_1(0) = 2p$ and $V_1(\pm 1) = p(1+k) + (1-2p) = pk + 1 - p$. Now, if $k <
(3p -1)/p$, then $V_1(-1) < V_1(0) > V_1(1)$ and hence the value function is
not quasi-convex. Hence, condition (M5) is necessary.

\section{Conclusion}
In this paper we consider a Markov decision process with continuous or
discrete state and action spaces and analyze the monotonicity of the optimal
solutions. In particular, we identify sufficient conditions under which the
value function and the optimal strategy are even and quasi-convex. The proof
relies on a folded representation of the Markov decision process and uses
stochastic monotonicity and submodularity. We present an example of optimal
power allocation in remote estimation and show that the sufficient conditions
are easily verified.

Establishing that the value function and optimal strategy are even and
quasi-convex has two benefits. First, such structured strategies are easier to
implement. Second, the structure of the value function and optimal strategy may
be exploited to efficiently solve the dynamic program.

For example, when the action space is discrete, say $|\ALPHABET U|=m$, then
even and quasi-convex strategy is characterized by $m-1$ thresholds. Such a
threshold-based strategy is simpler to implement than an arbitrary strategy.
Furthermore, the threshold structure also simplifies the search of the
optimal strategy. For discrete state spaces, see the monotone dynamic
programming presented in Sec.~\ref{sec:monotoneDP}; for continuous state
spaces, see~\cite{JC-JS-AM-ACC17}, where a simulation based algorithm is
presented to compute the optimal thresholds in remote estimation over a
packet drop channel.

Even for continuous action spaces, it is easier to search within the class of
even and quasi-convex strategies. Typically, some form of approximation is
needed to search for an optimal strategy. Two commonly used approximation
schemes are discretizing the action space or projecting the policy on to a
parametric family of function. If the action state is discretized, then the
search methods for discrete action spaces may be used. If the strategy is
projected on to a parametric family of function, then the structure may help
in reducing the size of the parameter space. For example, when approximating
an even and quasi-convex policy as a finite order polynomial, one can
restrict attention to polynomials where the coefficients of even powers are
positive and the coefficients of odd powers are zero.

In this paper, we assumed that the state space $\ALPHABET X$ was a subset of
reals. It will be useful to generalize these results to higher dimensions. 

\appendices

\section{Proof of Claim~\ref{clm:prop}}\label{app:proof_clm_prop}

We first prove some intermediate results:

\begin{lemma} \label{lem:basic}
  Under \textup{(M4)} and \textup{(M5)}, for any $x,y \in \ALPHABET X_{\ge
  0}$, we have that
  \[
    \varphi(y - x) \ge \varphi (y + x)
  \]
\end{lemma}
\begin{proof}
  We consider two cases: $y \ge x$ and $y < x$. 
  \begin{enumerate}
    \item If $y \ge x$, then $y + x \ge y - x \ge 0$. Thus, (M5) implies
      that $\varphi(y + x) \ge \varphi(y - x)$.
    \item If $y < x$, then $y + x \ge x - y$. Thus, (M5) implies that
      $\varphi(y + x) \ge \varphi(x - y) = \varphi(y - x)$, where the last
      equality follows from (M4).
  \end{enumerate}
\end{proof}

Some immediate implications of Lemma~\ref{lem:basic} are the following.
\begin{lemma} \label{lem:mono}
  Under \textup{(M4)} and \textup{(M5)}, for any $a \in \ALPHABET X$ and $x,y
  \in \ALPHABET X_{\ge 0}$, we have that
  \[
    a \big[ \varphi(y - ax) - \varphi (y + ax) \big] \ge 0.
  \]
\end{lemma}
\begin{proof}
  For $a \in \ALPHABET X_{\ge 0}$, from Lemma~\ref{lem:basic} we get that
  $\varphi(y - ax) \ge \varphi(y + ax)$. 
  For $a \in \ALPHABET X_{< 0}$, from Lemma~\ref{lem:basic} we get that
  $\varphi(y + ax) \ge \varphi(y - ax)$.
\end{proof}

\begin{lemma} \label{lem:mono2}
  Under \textup{(M4)} and \textup{(M5)}, for any $a, b, x,y
  \in \ALPHABET X_{\ge 0}$, we have that
  \[
    \varphi(y - ax - b) \ge \varphi (y + ax + b)
    \ge \varphi(y + ax + b + 1).
  \]
\end{lemma}
\begin{proof}
  By taking $y = y - b$ and $x = ax$ in Lemma~\ref{lem:basic}, we get
  \[
    \varphi(y - b - ax) \ge \varphi(y - b + ax).
  \]
  Now, by taking $y = y + ax$ and $x = b$ in Lemma~\ref{lem:basic}, we get
  \[
    \varphi(y + ax - b) \ge \varphi(y + ax + b).
  \]
  By combining these two inequalities, we get
  \[
    \varphi(y - ax - b) \ge \varphi(y + ax + b).
  \]
  The last inequality in the result follows from (M5).
\end{proof}

\begin{lemma}\label{lem:mono3}
  Under \textup{(M4)} and \textup{(M5)}, for $a \in \mathbb Z$ and $x,y \in
  \mathbb Z_{\ge 0}$, 
  \[ 
    \Phi (y+ax) + \Phi (y-ax) \ge \Phi (y + ax+ a) + \Phi (y -ax -a),
  \]
  where $\Phi$ is the cdf (cumulative distribution function) of $\varphi$.
\end{lemma}
\begin{proof}
  The statement holds trivially for $a = 0$. Furthermore, the statement does
  not depend on the sign of $a$. So, without loss of generality, we assume
  that $a > 0$. 

  Now consider the following series of inequalities (which follow from
  Lemma~\ref{lem:mono2})
  \begin{align*}
    \varphi(y - ax) &\ge \varphi(y + ax + 1),
    \\
    \varphi(y - ax - 1) &\ge \varphi(y + ax + 2),
    \\
    \cdots &\ge \cdots
    \\
    \varphi(y - ax -a + 1) &\ge \varphi(y + ax + a).
  \end{align*}
  Adding these inequalities, we get
  \[
    \Phi(y - ax) - \Phi(y - ax - a) 
    \ge
    \Phi(y + ax + a) - \Phi(y + ax),
  \]
  which proves the result.
\end{proof}

\subsection*{Proof of Claim~\ref{clm:prop}}

First, let's assume that $\ALPHABET X = \reals$. We prove each part
separately. 
\begin{enumerate}
  \item Fix $u \in [0, u_{\max}]$. $c(\cdot, u)$ is even because $d(\cdot)$
    is even (from (M3)). $c(\cdot, u)$ is quasi-convex because $1 - q(u)
    \ge 0$ (from (M0)) and $d(\cdot)$ is quasi-convex (from (M3)).
  \item Consider $e_1, e_2 \in \reals_{\ge 0}$ and $u_1, u_2 \in [0,
    u_{\max}]$ such that $e_1 \ge e_2$ and $u_1 \ge u_2$. The per-step cost
    is submodular on $\reals_{\ge 0} \times [0, u_{\max}]$ because
    \begin{align*}
      c(e_1, u_2) - c(e_2, u_2) &= (1 - q(u_2))( d(e_1) - d(e_2) ) \\
      & \stackrel{(a)}\ge (1 - q(u_1)) (d(e_1) - d(e_2)) \\
      &= c(e_1, u_1) - c(e_2, u_1),
    \end{align*}
    where $(a)$ is true because $d(e_1) - d(e_2) \ge 0$ (from (M3)) and
    $1 - q(u_2) \ge 1 - q(u_1) \ge 0$ (from (M0) and (M2)).
  \item Fix $u \in [0, u_{\max}]$ and consider $e, e_{+} \in \reals$.
    Then, $p(u)$ is even because
    \begin{align*}
      p(-e_+| {-e}; u) &= q(u) \varphi(e_+) + (1 - q(u)) \varphi(-e_{+} +
      ae) \\
      &\stackrel{(b)}= q(u) \varphi(e_+) + (1 - q(u)) \varphi(e_+ - ae)
      \\
      &= p(e_+ | e; u),
    \end{align*}
    where $(b)$ is true because $\varphi$ is even (from (M4)).
  \item First note that 
    \begin{align*}
      \hskip 1em & \hskip -1em
      S(y|x;u) = 1 - \int_{-\infty}^y \big[
      p(z|x;u) + p(-z|x;u) \big] dz \\
      &= 1 - \int_{-\infty}^y q(u) \big[ \varphi(z) + \varphi(-z) \big]dz \\
      &\quad - \int_{-\infty}^y (1 - q(u)) 
      \big[ \varphi(z - ax) + \varphi(-z - ax) ] dz
      \\
      &\stackrel{(c)}=
      1 - 2 q(u) \Phi(y) \\
      &\quad - (1 - q(u))\big[
      \Phi(y - ax) + \Phi(y + ax) \big]
    \end{align*}
    where $\Phi$ is the cumulative distribution of $\varphi$ and $(c)$ uses the
    fact that $\varphi$ is even (condition (M4)).

    Let $S_x(y|x;u)$ denote $\partial S/\partial x$. Then
    \[
      S_x(y|x;u) = (1 - q(u)) a
      \big[ \varphi(y - ax) - \varphi(y + ax) \big].
    \]
    From (M0) and Lemma~\ref{lem:mono}, we get that $S_x(y|x;u) \ge 0$ for
    any $x,y \in \reals_{\ge 0}$ and $u \in [0, u_{\max}]$. Thus,
    $S(y|x;u)$ is increasing in $x$. 

    Furthermore, from (M2) $S_x(y|x;u)$ is decreasing in $u$. Thus,
    $S(y|x;u)$ is submodular in $(x,u)$ on $\reals_{\ge 0} \times [0,
    u_{\max}]$. 
\end{enumerate}

Now, let's assume that $\ALPHABET X = \mathbb Z$. The proof of the first
three parts remains the same. Now, in part 4), it is still the case that
\begin{multline*}
  S(y|x;u) = 1 - 2q(u) \Phi(y)\\
  - (1 - q(u))\big[\Phi(y - ax) + \Phi(y + ax) \big]
\end{multline*}
However, since $\ALPHABET X$ is discrete, we cannot take the partial
derivative with respect to $x$. Nonetheless, following the same intuition,
for any $x,y \in \mathbb Z_{\ge 0}$, consider
\begin{multline}
  S(y|x+1;u) - S(y|x;u) = (1-q(u)) \big[
    \Phi(y + ax) \\
    {}- \Phi(y + ax +a) + \Phi(y - ax) - \Phi(y - ax - a)
  \big]
  \label{eq:S-discrete}
\end{multline}

Now, by Lemma~\ref{lem:mono3}, the term in the square bracket is positive,
and hence $S(y|x;u)$ is increasing in $x$. Moreover, since $(1 - q(u))$ is
decreasing in $u$, so is $S(y|x+1;u) - S(x|x;u)$. Hence, $S(y|x;u)$ is
submodular in $\ALPHABET Z_{\ge 0} \times [0, u_{\max}]$.

\bibliographystyle{IEEEtran}
\bibliography{IEEEabrv,monotone_ref}

\end{document}